\documentclass[12pt]{article}

\newcommand{\JMPcite}[1]{\kern-.2em${}^{\hbox{\scriptsize\cite{#1}}}$}%

\newcommand{\be}{\begin{equation}}
\newcommand{\ee}{\end{equation}}
\newcommand{\bea}{\begin{eqnarray}}
\newcommand{\eea}{\end{eqnarray}}

\def\bbbc{{\mathchoice {\setbox0=\hbox{$\displaystyle\rm C$}\hbox{\hbox
to0pt{\kern0.4\wd0\vrule height0.9\ht0\hss}\box0}}
{\setbox0=\hbox{$\textstyle\rm C$}\hbox{\hbox
to0pt{\kern0.4\wd0\vrule height0.9\ht0\hss}\box0}}
{\setbox0=\hbox{$\scriptstyle\rm C$}\hbox{\hbox
to0pt{\kern0.4\wd0\vrule height0.9\ht0\hss}\box0}}
{\setbox0=\hbox{$\scriptscriptstyle\rm C$}\hbox{\hbox
to0pt{\kern0.4\wd0\vrule height0.9\ht0\hss}\box0}}}}
\def\bbbc{{\rm I\!C}}

\begin{document}

\centerline{\bf{Gibbs Measures For SOS Models On a Cayley Tree}}
\vskip 0.5 truecm
\centerline{\bf{U.A. Rozikov\footnote{Institute of Mathematics,
Uzbek Academy of Sciences, Tashkent 700143, Uzbekistan} and Y.M. Suhov
\footnote{Statistical Laboratory, DPMMS, University of Cambridge,
Cambridge CB3 0WB, UK}}}
\vskip 0.5 truecm

-----------------------------------------------

{\bf Abstract.} We consider a nearest-neighbor SOS (solid-on-solid)
model, with several spin values $0,1,\ldots ,m$, $m\geq 2$, and
zero external field, on a Cayley tree of order $k$ (with
$k+1$ neighbors). The SOS model can be treated as a
natural generalisation of the Ising model (obtained for
$m=1$). We mainly assume that $m=2$ (three spin values) and study
translation-invariant (TI) and `splitting'
(S) Gibbs measures (GMs). [Splitting GMs have
a particular Markov-type property specific for a tree.]
Furthermore, we focus on symmetric TISGMs, with respect
to a `mirror' reflection of the spins. [For the Ising
model (where $m=1$), such measures are reduced to the `disordered' phase
obtained for free boundary conditions, see [BRZ], [M1,2],
[MP].] For $m=2$, in the anti-ferromagnetic (AFM) case, a  
symmetric TISGM (and
even a general TISGM) is unique for all
temperatures. In the ferromagnetic (FM) case, for $m=2$, the number of
symmetric TISGMs and (and the number of general TISGMs) 
varies with the temperature: this gives an
interesting example of phase transition. Here we identify
a critical inverse temperature,
$\beta^1_{\rm{cr}}$ ($=T_{\rm{cr}}^{\rm{STISG}}$)
$\in (0,\infty )$ such that $\forall$
$0\leq \beta\leq\beta^1_{\rm{cr}}$, there
exists a unique symmetric TISGM $\mu^*$ and $\forall$
$\beta >\beta^1_{\rm{cr}}$ there
are exactly three symmetric TISGMs: $\mu^*_+$ (a `bottom'
symmetric TISGM), $\mu^*_{\rm m}$ (a `middle' symmetric
TISGM) and $\mu^*_-$
(a `top' symmetric TISGM).
For $\beta>\beta^1_{\rm{cr}}$ we also
construct a continuum of distinct, symmertric SGMs which are non-TI.

Our second result gives complete description of the set of periodic
Gibbs measures for the SOS model on a Cayley tree. A complete description
of periodic GMs means a characterisation of such measures with 
respect to any given normal subgroup of finite index in the representation 
group of the tree. We show that (i) for an FM SOS model, for any 
normal subgroup of finite index, each periodic  SGM is
in fact TI. Further, (ii)   
for an AFM SOS model, for 
any normal subgroup of finite index, each periodic 
SGM is either TI or has 
period two (i.e., is a chess-board SGM).

-----------------------------------------------


\section{Introduction}
One of the central problems in the theory of Gibbs measures
(GMs)
is to describe infinite-volume (or limiting) GMs
corresponding to a given Hamiltonian. The existence of such
measures for a wide class of Hamiltonians was established
in the ground-breaking work of Dobrushin (see, e.g., Ref. [1]).
However, a complete analysis of the set of limiting GMs
for a specific Hamiltonian is often a difficult
problem. On a cubic lattice, for small values of $\beta={1\over T}$,
where $T>0$ is the temperature, a GM is unique
(Refs [1-3]) which reflects a physical fact that at high temperatures
there is no phase transitions. The analysis for
low temperatures requires specific assumptions on
the form of the Hamiltonian.

In this paper we consider models with a nearest neighbour
interaction on a Cayley tree (CT). Models on a CT
were discussed in Refs. [2], [4]--[6]. A classical example of
such a model is the Ising model, with two values of spin, $\pm 1$.
It was considered in Refs. [2], [6] and became a focus of active
research in the first half of the 1990's and afterwards; see Refs [7]--[13].
Models considered in the present paper are generalisations
of the Ising model and can be described as SOS (solid-on-solid)
models with constraints; see below. In the case of a cubic
lattice they were analysed in Ref. [14] where an analogue of
the so-called Dinaburg--Mazel--Sinai theory was developed.
Besides interesting phase transitions in these models, the
attention to them is motivated by applications, in particular
in the theory of communication networks; see, e.g., Refs [15].

A CT ${\cal T}^k=(V,A)$ of order $k\geq 1$
is an infinite homogeneous tree, i.e., a graph without
cycles, with exactly $k+1$ edges incident to each vertex.
Here $V$ is the set of vertices and $A$ that of edges
(arcs).

We consider models where the spin takes values in the set
$\Phi:=\{0,1,\ldots , m\}$, $m\geq 2$, and is assigned to the vertices
of the tree. A configuration $\sigma$ on $V$ is then defined
as a function $x\in V\mapsto\sigma (x)\in\Phi$;
the set of all configurations is $\Phi^V$.
The (formal) Hamiltonian is of an SOS form:
$$ H(\sigma)=-J\sum_{\langle x,y\rangle\in L}
|\sigma(x)-\sigma(y)|,\eqno (1.1)$$
where $J\in R$ is a coupling constant. As usually,
$\langle x,y\rangle$ stands for nearest neighbor vertices.

The SOS model of this type can be considered as a generalisation
of the Ising model (which arises when $m=1$). Here, $J<0$
gives a ferromagnetic (FM) and $J>0$ an antiferromagnetic (AFM)
model. In the FM case the ground states are `flat' configurations,
with $\sigma (x)\equiv j\in\Phi$ (there are $m+1$ of them),
in the AFM two `contrasting'
checker-board configurations where $|\sigma (x)-\sigma (y)|=m$
$\forall$ $\langle x,y\rangle$. Compared with
the Potts model (see, e.g., [16]--[19]), the SOS has
`less symmetry' and therefore more diverse structure
of phases. For example, in the FM case
it is intuitively plausible that the ground states
corresponding to `middle-level'
surfaces will be `dominant'. This observation was made formal
in [14] for the model on a cubic lattice.

We consider a standard sigma-algebra ${\cal B}$ of subsets of
$\Phi^V$ generated by cylinder subsets; all probability
measures are considered on $(\Phi^V,{\cal B})$.
A probability measure $\mu$ is called a GM (with
Hamiltonian $H$) if it satisfies the DLR equation: $\forall$
$n=1,2,\ldots$ and $\sigma_n\in\Phi^{V_n}$:
$$\mu\left(\left\{\sigma\in\Phi^V :\;
\sigma\big|_{V_n}=\sigma_n\right\}\right)=
\int_{\Phi^V}\mu ({\rm d}\omega)\nu^{V_n}_{\omega|_{W_{n+1}}}
(\sigma_n),\eqno (1.2)$$
where $\nu^{V_n}_{\omega|_{W_{n+1}}}$ is the conditional
probability:
$$ \nu^{V_n}_{\omega|_{W_{n+1}}}(\sigma_n)=\frac{1}{Z_n\left(
\omega\big|_{W_{n+1}}\right)}\exp\;\left(-\beta H
\left(\sigma_n\,||\,\omega\big|_{W_{n+1}}\right)\right).
\eqno (1.3)$$
Here and below, $W_l$ stands for a `sphere' and $V_l$ for a
`ball' on the tree, of radius $l=1,2,\ldots$,
centered at a fixed vertex $x^0$ (an origin):
$$W_l=\{x\in V: d(x,x^0)=l\},\;\;V_l=\{x\in V: d(x,x^0)\leq l\};$$
distance $d(x,y)$, $x,y\in V$, is the length of
(i.e. the number of edges in)
the shortest path connecting $x$ with $y$. $\Phi^{V_n}$
is the set of configurations
in $V_n$ (and $\Phi^{W_n}$ that in $W_n$; see below). Furthermore,
$\sigma\big|_{V_n}$ and
$\omega\big|_{W_{n+1}}$ denote the restrictions
of configurations $\sigma,\omega\in\Phi^V$ to $V_n$
and $W_{n+1}$, respectively. Next,
$\sigma_n:\;x\in V_n\mapsto \sigma_n(x)$ is a configuration in $V_n$
and
$H\left(\sigma_n\,||\,\omega\big|_{W_{n+1}}\right)$
is defined as the sum $H\left(\sigma_n\right)+U\left(\sigma_n,
\omega\big|_{W_{n+1}}\right)$ where
$$H\left(\sigma_n\right)
=-J\sum_{\langle x,y\rangle\in L_n}|\sigma_n(x)-\sigma_n(y)|,\;
U\left(\sigma_n,
\omega\big|_{W_{n+1}}\right)=
-J\sum_{\langle x,y\rangle:\;x\in V_n,y\in W_{n+1}}
|\sigma_n(x)-\omega (y)|,\eqno (1.4)$$
and
$$L_n=\{\langle x,y\rangle\in L: x,y\in V_n\}.$$
Finally, $Z_n\left(\omega\big|_{W_{n+1}}\right)$
stands for the partition function in $V_n$, with
the boundary condition $\omega\big|_{W_{n+1}}$:
$$Z_n\left(\omega\big|_{W_{n+1}}\right)=
\sum_{{\widetilde\sigma}_n\in\Phi^{V_n}}
\exp\;\left(-\beta H
\left({\widetilde\sigma}_n\,||\,\omega
\big|_{W_{n+1}}\right)\right).\eqno (1.5)$$

Because of the nearest-neghbour character of the interaction,
the GMs possess a natural Markov property: given
a configuration $\omega_n$ on $W_n$, random configurations
in $V_{n-1}$ (i.e., `inside' $W_n$) and in $V\setminus V_{n+1}$
(i.e., `outside' $W_n$) are conditionally independent.
It is known (see, e.g., [1], [2]) that $\forall$
$\beta >0$, the GMs form a non-empty
convex compact set in the space of probability measures.
Extreme measures, i.e., extreme points of this set are
associated with pure phases.
Furthermore, any GM is an integral of
extreme ones (the extreme decomposition).
It is true that for any sequence of configurations
$\omega^{(n)}\in\Phi^V$,
every limiting point of measures ${\widetilde\nu}^{V_n}_{
\omega^{(n)}|_{W_{n+1}}}$ is a GM. Here, for a given
$\omega\in\Phi^V$,
${\widetilde\nu}^{V_n}_{\omega|_{W_{n+1}}}$
is a measure on $\Phi^V$ such that $\forall$ $n'>n$:
$${\widetilde\nu}^{V_n}_{\omega|_{W_{n+1}}}\left(\left\{
\sigma\in\Phi^V:\;\sigma\big|_{V_{n'}}=
\sigma_{n'}\right\}\right)=\left\{\begin{array}{ll}
\displaystyle{\nu^{V_n}_{\omega|_{W_{n+1}}}\left(\sigma_{n'}
\big|_{V_n}\right),}&\hbox{if }\;\sigma_{n'}\big|_{
V_{n'}\setminus V_n}=\omega\big|_{V_{n'}\setminus V_n},\\
0,&\hbox{otherwise.}\end{array}\right.\eqno (1.6)$$
The converse is also true: every GM $\mu$ can be obtained
as a limiting point for measures
${\widetilde\nu}^{V_n}_{\omega^{(n)}|_{W_{n+1}}}$ with a suitable
sequence of configurations $\omega^{(n)}\Phi^V$. We call
such a sequence $\omega^{(n)}$ the boundary conditions for
GM $\mu$.

We use a standard definition of a translation-invariant (TI)
measure (see, e.g., [4]). Also, call measure $\mu$ symmetric
(S) if it is preserved under the simultaneous change $j\mapsto m-j$
at each vertex $x\in V$. The main object of study in this
paper are symmetric TI mesaures.

An important role is played by a specific monotonicity
displayed by the FM model (with $J<0$). Namely, write
$\sigma\leq\sigma'$ if configurations $\sigma$ and
$\sigma'$ obey $\sigma (x)\leq\sigma'(x)$ $\forall$
$x\in V$. This partial order defines a concept of
a monotone increasing and monotone decreasing function
$f:$ $\Phi^V\to{\mathbf R}$. For two probability measures
$\mu_1$ and $\mu_2$ we then write $\mu_1\leq\mu_2$
if $\int fd\mu_1\leq\int fd\mu_2$ for each monotone increasing
$f$. It turns out that for the `extreme' configurations,
$\omega^0$ with $\omega^0(x)\equiv 0$ and $\omega^2$ with
$\omega^2(x)\equiv 2$, there exist the limits
$\nu^0=\lim_{n\to\infty}{\widetilde\nu}^{V_n}_{\omega^0|_{W_{n+1}}}$
and $\nu^2=\lim_{n\to\infty}{\widetilde\nu}^{V_n}_{\omega^2|_{W_{n+1}}}$
(both measure sequences are monotone).
$\nu^0$, $\nu^2$ are TIGMs and possess the following
minimality and maximality properties: $\nu^1\leq\mu\leq \nu^2$ $\forall$
GM $\mu$. Because of that, they are both extreme (although
not symmetric).
The question of whether a GM is non-unique is then reduced to
whether $\nu^1=\nu^2$. However, finer properties of GMs
require further specifications.

\section{Construction of splitting GMs}

Following Ref. [2] (and subsequent papers [5-13]), we consider
a special class of GMs. These measures are called
in Ref. [2] Markov chains and in Refs [5], [6] entrance laws.
In this paper we call them splitting GMs, to emphasize
the property that, in addition to the aforementioned Markov property,
they satisfy the following condition: given a configuration $\sigma_n$
in $V_n$, the values $\sigma (y)$ at sites $y\in W_{n+1}$
are conditionally independent.

Write $x<y$ if the path from $x^0$ to $y$ goes
through $x$. Call vertex $y$ a direct successor
of $x$ if $y>x$ and $x,y$ are nearest neighbours.
Denote by $S(x)$ the set of direct successors of $x$.
Observe that any vertex $x\ne x^0$ has $k$ direct
successors and $x^0$ has $k+1$.

Let $h:\;x\mapsto h_x=(h_{0,x}, h_{1,x},...,h_{m,x})
\in R^{m+1}$ be a real vector-valued function of $x\in V\setminus
\{x^0\}$. Given $n=1,2,\ldots$,
consider the probability distribution $\mu_n$ on
$\Phi^{V_n}$ defined by
$$\mu^{(n)}(\sigma_n)=Z_n^{-1}\exp\left(-\beta H(\sigma_n)
+\sum_{x\in W_n}h_{\sigma(x),x}\right), \eqno (2.1)$$
Here, as before, $\sigma_n:x\in V_n\mapsto \sigma(x)$
and $Z_n$ is the corresponding partition function:
$$Z_n=\sum_{{\widetilde\sigma}_n\in\Phi^{V_n}}
\exp\left(-\beta H({\widetilde\sigma}_n)
+\sum_{x\in W_n}h_{{\widetilde\sigma}(x),x}\right).\eqno (2.2)$$

We say that the probability distributions $\mu^{(n)}$
are compatible if $\forall$ $n\geq 1$ and $\sigma_{n-1}\in\Phi^{V_{n-1}}$:
$$\sum_{\omega_n\in\Phi^{W_n}}\mu^{(n)}(\sigma_{n-1}\vee\omega_n)=
\mu^{(n-1)}(\sigma_{n-1}).\eqno (2.3)$$
Here $\sigma_{n-1}\vee\omega_n\in\Phi^{V_n}$ is the concatenation
of $\sigma_{n-1}$ and $\omega_n$.
In this case there exists a unique measure $\mu$ on
$\Phi^V$ such that, $\forall$ $n$ and
$\sigma_n\in\Phi^{V_n}$, $\mu \left(\left\{\sigma
\Big|_{V_n}=\sigma_n\right\}\right)=\mu^{(n)}(\sigma_n)$. Such
a measure is
called a splitting GM (SGM) corresponding to
Hamiltonian $H$ and function $x\mapsto h_x$, $x\neq x^0$.

The following statement describes conditions on
$h_x$ guaranteeing
compatibility of distributions $\mu^{(n)}(\sigma_n).$
\vskip 0.5 truecm

{\bf Proposition 2.1.} {\sl Probability distributions
$\mu^{(n)}(\sigma_n)$, $n=1,2,\ldots$, in}
(2.1) {\sl are compatible iff for any $x\in V\setminus\{x^0\}$
the following equation holds:
$$ h^*_x=\sum_{y\in S(x)}F(h^*_y,m,\theta). \eqno(2.4)$$
Here, and below
$$\theta=\exp(J\beta ),\eqno (2.5)$$
$h^*_x$ is stands for the vector
$(h_{0,x}-h_{m,x}, h_{1,x}-h_{m,x},...,h_{m-1,x}-h_{m,x})$ and
the vector function $F(\;\cdot\;,m,\theta ):\;R^m\to R^m$
is $F(h,m,\theta )=
(F_0(h,m,\theta ),\ldots ,F_{m-1}(h,m,\theta))$, with}
$$F_i(h,m,\theta )=\ln{\sum_{j=0}^{m-1}
\theta^{|i-j|}\exp(h_j)+\theta^{m-i}\over
\sum_{j=0}^{m-1}\theta^{m-j}\exp(h_j)+1},\;h=(h_0,h_1,...,h_{m-1}),
i=0,\ldots ,m-1. \eqno (2.6)$$
\vskip 0.5 truecm

{\bf Proof.} {\sl Necessity} (cf. [16]). Suppose that (2.3)
holds; we want to prove (2.4). Substituting (2.1) in (2.3),
obtain, $\forall$ configurations $\sigma_{n-1}$:
$x\in V_{n-1}\mapsto\sigma_{n-1}(x)\in\Phi$:
$$\frac{Z_{n-1}}{Z_n}\sum_{\omega_n\in\Phi^{W_n}}
\exp\left(\sum_{x\in W_{n-1}}\sum_{y\in S(x)}
(J\beta |\sigma_{n-1}(x)-\omega_n(y)|+h_{\omega_n(y),y})\right)=
\exp\left(\sum_{x\in W_{n-1}}h_{\sigma_{n-1}(x),x}\right),
\eqno (2.7)$$
where $\omega_n$: $x\in W_n\mapsto\omega_n(x)$.

From (2.7) we get:
$${Z_{n-1}\over Z_n}\sum_{\omega_n\in\Phi^{W_n}}
\prod_{x\in W_{n-1}}\prod_{y\in S(x)}
\exp\,(J\beta |\sigma_{n-1}(x)-\omega_n(y)|+
h_{\omega_n(y),y})=\prod_{x\in W_{n-1}}
\exp\,(h_{\sigma_{n-1}(x),x}). \eqno (2.8)$$
Consequently, $\forall$ $i\in\Phi$,
$$\prod_{y\in S(x)}\frac{\sum_{j\in\Phi}\exp\,(J\beta |i-j|+
h_{j,y})}{\sum_{j\in\Phi}\exp\,(J\beta |m-j|+h_{j,y})}=
\exp\,(h_{i,x}-h_{m,x}). \eqno (2.9)$$
Introducing $\theta$ as in (2.5) and denoting
$h^*_{i,x}=h_{i,x}-h_{m,x}$, we get (2.4) from (2.9).

{\sl Sufficiency.} From (2.4) we obtain
(2.9), (2.8) and  (2.7) i.e. (2.3). The proof is complete.
\vskip 0.5 truecm

{\bf Proposition 2.2.} {\sl Any measure $\mu$ with local
distributions $\mu^{(n)}$ satisfying} (2.1), (2.3) {\sl is an
SGM.}
\vskip 0.5 truecm

{\bf Proof.} Straightforward.
\vskip 0.5 truecm

{\bf Proposition 2.3.} {\sl An SGM $\mu$ is TI iff
$h_{j,x}$ does not depend on $x$: $h_{j,x}\equiv h_j$, $x\in V$,
$j\in\Phi$, and symmetric TI iff $h_j=h_{m-j}$, $j\in\Phi$.}
\vskip 0.5 truecm

{\bf Proof.} Straightforward.
\vskip 0.5 truecm

{\bf Proposition 2.4.} {\it Any extreme GM is an SGM.}
\vskip 0.5 truecm

{\bf Proof.} See Ref [4], Theorem 12.6.

\section{The critical value $\beta^1_{\rm{cr}}$}

From Proposition 2.2 it follows that for any $h=\{h_x,\ \ x\in V\}$
satisfying (2.4) there exists a unique GM $\mu$ (with
restrictions $\mu^{(n)}$ as in (2.1)) and vice versa. However,
the analysis of solutions to (2.4) for an arbitrary
$m$ is not easy. We now
suppose that the number of spin values $m+1$ is 3 i.e. $m=2$
and $\Phi =\{0,1,2\}$. Throughout the paper we assume that
$h_{2,x}\equiv 0$ ($h_{m,x}\equiv 0$ for general $m$).

It is natural to begin with translation-invariant solutions
where $h_x=h\in R^m$ is constant. Unless stated otherwise,
we concentrate on the simplest
case where $m=2$, i.e. spin values are $0$, $1$ and $2$.
In this case we obtain from (2.4), (2.5):
$$h_{0,x}=\sum_{y\in S(x)}\ln{\exp(h_{0,y})+
\theta\exp(h_{1,y})+\theta^2\over
\theta^2\exp(h_{0,y})+\theta\exp(h_{1,y})+1},\;\;
h_{1,x}=\sum_{y\in S(x)}\ln{\theta\exp(h_{0,y})+\exp(h_{1,y})+\theta
\over
\theta^2\exp(h_{0,y})+\theta\exp(h_{1,y})+1}.\eqno (3.1)$$

Set $z_0=\exp(h_{0,x})$, $z_1=\exp(h_{1,x})$ (and $z_2=1$), $x\in V$.
From (3.1) we have
$$ z_0=\left({z_0+\theta z_1+\theta^2 \over \theta^2z_0+\theta z_1+1}
\right)^k,\eqno (3.2.a)$$
$$ z_1=\left({\theta z_0+z_1+\theta \over \theta^2z_0+\theta z_1+1}
\right)^k.\eqno (3.2.b)$$
Observe that $z_0=1$ satisfies equation
(3.2.a) independently of
$k$, $\theta$ and $z_1$.
Substituting $z_0=1$ into (3.2.b), we obtain
$$ z_1=\left({2\theta +z_1 \over \theta^2+\theta z_1+1}\right)^k.
\eqno (3.3)$$

Set:
$$a=2\theta^{k+1},\;b = {1+\theta^2\over 2\theta^2},
\;\;x={z_1\over 2\theta}.\eqno (3.4)$$
Then from (3.3):
$$ ax=\left({1+x \over b+x}\right)^k. \eqno (3.5)$$
In Proposition 3.1 below we analyse solutions to equation (3.5)
with independently varying parameters $a,b>0$. The
proof of Proposition 3.1
repeats an argument from [2], Proposition 10.7.
\vskip 0.5 truecm

{\bf Proposition 3.1.} {\sl Equation} (3.5) {\sl with $x\geq 0$,
$k\geq 1$, $a,b >0$ has a unique solution if either $k=1$ or
$b\leq ({k+1 \over k-1})^2$. If $k>1$ and $b>({k+1 \over k-1})^2$
then there exist $\nu_1(b,k)$, $\nu_2(b,k)$, with
$0<\nu_1(b,k)< \nu_2(b,k)$, such that the equation has three
solutions if $\nu_1(b,k)<a< \nu_2(b,k)$ and has two if either
$a=\nu_1(b,k)$ or $a= \nu_2(b,k)$. In fact:
$$ \nu_i(b,k)={1\over x_i}\left({1+x_i \over b+x_i}\right)^k, $$
where $x_1,x_2$ are the solutions of}
$$ x^2+[2-(b-1)(k-1)]x+b=0.$$
\vskip 0.5 truecm

Now consider $a$ and $b$ as functions of $\beta$ (for a fixed $J$
as specified in (3.4) and (2.6)).
\vskip 0.5 truecm

{\bf Proposition 3.2.} {\sl If $J\geq 0$ then the system of
equations} (3.2.a), (3.2.b) {\sl has a unique solution.}
\vskip 0.5 truecm

{\bf Proof.} Let $A=z_0+\theta z_1+\theta^2$, $B=\theta^2z_0+\theta z_1+1$,
then from (3.2.a) we have:
$$ (z_0-1)[B^k+(\theta^2-1)(A^{k-1}+...+B^{k-1})]=0 \eqno (3.6)$$
Since $\theta\geq 1$ ($J\geq 0$), we deduce from (3.6) that $z_0=1$
is the only solution. Then
$b={1+\theta^2\over 2\theta^2}\leq 1<({k+1\over k-1})^2$.
By Proposition 3.1, equation (3.3) has a unique solution. Thus we have
proved that system (3.2.a), (3.2.b) has a unique solution.
\vskip 0.5 truecm

{\bf Proposition 3.3.} {\sl If $J<0$ then for $\beta\leq\frac{1}{2J}
\ln\frac{(k-1)^2}{k^2+6k+1}$, the system of equations} (3.2.a),
(3.2.b) {\sl has a unique solution of the form $(1,z^*)$ (i.e., a unique
solution $(z^*_0,z^*_1)$ with $z^*_0=1$) and for $\beta >
\frac{1}{2J}\ln\frac{(k-1)^2}{k^2+6k+1}$, precisely
three such solutions, $(1, z^*_{1,-})$,
$(1, z^*_{1,{\rm m}})$,
$(1, z^*_{1,+})$, with  $0<z^*_{1,-}<z^*_{1,{\rm m}} <z^*_{1,+}$ and
$z^*_{1,i}=e^{h^*_{1,i}}$, $i=-,{\rm m},+$, (see} (3.1) {\sl ).}
\vskip 0.5 truecm

{\bf Proof.} The value $\frac{1}{2J}\ln\frac{(k-1)^2}{
k^2+6k+1}$ is the solution of equation $b=({k+1\over k-1})^2$.
Other statements of Proposition 3.3 are consequences of Proposition
3.1.
\vskip 0.5 truecm

For brevity we say that $z^*$ and $z^*_-$, $z^*_{\rm m}$, $z^*_+$
give symmetric solutions to (3.2.a,b).
\vskip 0.5 truecm

{\bf Definition.} In the FM case, set:
$$\beta^1_{\rm{cr}}=\frac{1}{2J}
\ln\frac{(k-1)^2}{k^2+6k+1}>0.\eqno (3.7)$$
\vskip 0.5 truecm

Going back to (2.5), summarise:
\vskip 0.5 truecm

{\bf Theorem 1.} {\sl For the AFM SOS model, with $J>0$
and $m=2$, the TISGM exists and is unique $\forall$ $\beta\geq
0$. In fact, it is a symmetric TISGM.}

{\sl For the FM SOS model, with $J<0$
and $m=2$:}

1) {\sl If $k\geq 2$ and $0\leq\beta\leq\beta^1_{\rm{cr}}$
then there exists a unique symmetric TISGM, $\mu^*$.}

2) {\sl If $k\geq 2$ and $\beta >\beta^1_{\rm{cr}}$
then there exist precisely three symmetric TISGMs
$\mu^*_-$, $\mu^*_{\rm m}$, $\mu^*_+$ corresponding to $h^*_i
=\ln\,z^*_i$, $i=-,{\rm m},+$.}
\vskip 0.5 truecm

{\bf Remark 1.} In the AFM case, the phase transition is
manifested in the break of the TI property. More precisely,
it is expected that for $\beta$ small there exists a
unique translation-periodic SGM (which is TI) while for $\beta$ large
there are several such measures. 
\vskip 0.5 truecm

In the FM case, observe that values  $z^*_{i}$,
$i=-,{\rm m},+$, vary
with $\beta$. It is easy to show that as $\beta\to\infty$,
$z^*_-\to 0$,
$z^*_{\rm m}\to 1$ and $z^*_+\to\infty$. Correspondingly,
we make a
\vskip 0.5 truecm

{\bf Conjecture 1.} For $m=2$, $k\geq 2$ and $J<0$,
as $\beta\to\infty$, measure
$\mu^*_-$ tends to the half-sum
$\frac{1}{2}\big(\delta_{\omega^0}+\delta_{\omega^2}\big)$,
$\mu^*_{\rm m}$ to
the mean $\frac{1}{3}\big(\delta_{\omega^0}+\delta_{\omega^1}+
\delta_{\omega^2}\big)$ and $\mu^*_+$
to $\delta_{\omega^1}$. Here $\delta_\omega$ stands for the Dirac
delta-measure sitting on configuration $\omega\in\Phi^V$
and $\omega^i$ has $\omega^i(x)\equiv i$, $i=0,1,2$.
\vskip 0.5 truecm

On the other hand, we can say that for $\beta\leq\beta^1_{\rm{cr}}$,
all three measures coincide and in the limit $\beta\to 0$ give
a Bernoulli measure,
with iid and equiprobable values $\sigma (x)=0,1,2$, $x\in V$.
\vskip 0.5 truecm

{\bf Remark 2.} Note that $\beta^1_{\rm{cr}}$ may
not be the first critical value of the inverse temperature
for the FM model.
Namely, there exists $\beta^0_{\rm{cr}}$
($=\beta^{\rm{TIGM}}_{\rm c}$) $\in
(0,\beta^1_{\rm{cr}}]$ such that (i)
for $0\leq \beta\leq\beta^0_{\rm{cr}}$, a
minimal GM, $\mu_-$, and a maximal, $\mu_+$,
coincide, and the whole set of GMs is reduced to 
a unique measure which is therefore extreme (and
coincides with symmetric TISGM $\mu^*$), (ii)
for $\beta >\beta^0_{\rm{cr}}$, $\mu_-$ and $\mu_+$
are distinct (they are always extreme TISGMs,
but not symmetric). Thus,
for $\beta\geq\beta^1_{\rm{cr}}$, there are five
TISGMs (in a natural order: $\mu_-\leq
\mu^*_-\leq\mu^*_{\rm m}\leq\mu^*_+\leq\mu_+$)
three of which are symmetric. 
It is not known whether $\beta^0_{\rm{cr}}
=\beta^1_{\rm{cr}}$ (it is our {\bf Conjecture 2}).
\vskip 0.5 truecm

The following Proposition 3.5 describes a useful property
of general (non-TI) solutions $h_x=(h_{0,x};h_{1,x})$ to (3.1)
with $h_{0,x}\equiv 0$ (or $z^*_0\equiv 1$). As before,
$h_{0,x}$ gives a solution to the first equation in (3.1),
regardless of $h_{1,x}$ and $\theta$.
\vskip 0.5 truecm

{\bf Proposition 3.5.} {\sl For $J<0$, $k\geq 2$
and $\beta >\beta^1_{\rm{cr}}$, if
$h_x=(0;h_{1,x})$ is a solution of} (3.1) {\sl then, with
$h_{1,x}=\ln z_{1,x}$,
$$z^*_{-}\leq z_{1,x}\leq z^*_{+},\ \ x\in V \eqno (3.8)$$
where $z^*_{-}< z^*_{+}$ are the symmetric solutions of}
(3.2.a,b) {\sl i.e., the solutions of} (3.3).
\vskip 0.5 truecm

{\bf Proof.} Denote $z_x=\exp(h_{1,x}).$ Then from (3.1) we get
$$ z_x=\prod^k_{i=1}{2\theta+z_{x_i}\over 1+\theta^2+\theta z_{x_i}}, \ \
z_{x_j}>0,\ \ j=1,\ldots, k,$$
where $x_j, \ \ j=1,\ldots ,k$ are direct successors of $x$.
Denote $\varphi(x,\theta)={2\theta +x\over 1+\theta^2+\theta x}.$
Consider
$$G(x_1,...,x_k)=\prod_{i=1}^k\varphi(x_i,\theta)
, \ \ x_i>0, \ \
i=1,\ldots ,k.$$
Set the map $x\mapsto \psi(x,\theta,k)=(\varphi(x,\theta))^k$.
Clearly, $\psi(0,\theta,k)\leq G(x_1,...,x_k)\leq
\psi(\infty,\theta,k)$. Thus for
$z_x$ we get
$\psi(0,\theta,k)\leq z_x\leq \psi(\infty, \theta, k)$.
Now consider $G(x_1,...,x_k)$ with
$\psi(0,\theta,k)\leq x_j\leq \psi(\infty,\theta, k)$.
Here we have
$$\psi(\psi(0,\theta, k),\theta, k)\leq z_x
\leq \psi(\psi(\infty,\theta, k),\theta, k).$$
Repeating this argument, we see that for the $n$th iteration
$\psi^{(n)}$ of $\psi$:
$$\psi^{(n)}(0,\theta,k)\leq z_x\leq \psi^{(n)}(\infty,\theta, k),$$
for all $n\geq 1$ and $x\in V\setminus\{x^0\}$. The sequence
$\psi^{(n)}(\infty,\theta, k)$ is
decreasing and bounded from below by $z^*_{+}$. Its limit is a
fixed point for $\psi$ and thus equal to $z^*_+$. The lower
bound for $z_x$ is similar and gives $z^*_{-}$.
\vskip 0.5 truecm

{\bf Proposition 3.6.} {\sl For $J<0$ and $\beta\leq
\beta^1_{\rm{cr}}$, measure $\mu^*$ is the only splitting
GM such that $z_{0,x}=0$, $x\in V\setminus\{x^0\}$ (regardless
whether it is TI or not). Thus, $\mu^*$ is the only symmetric
SGM.}
\vskip 0.5 truecm

{\bf Proof.} In this case equation (3.1) with $h_{0,x}=0$ has
a unique solution $h_x=(0,\ln z^*)$.
\vskip 0.5 truecm

{\bf Conjecture 3.} In the case $m=2$, $J<0$ and
$\beta\leq\beta^1_{\rm{cr}}$, $\mu^*$ is the unique GM
and hence extreme.
\vskip 0.5 truecm

{\bf Conjecture 4.} In the case $k\geq 2$, $J<0$ and
$\beta >\beta^1_{\rm{cr}}$,
the boundary condition for the top symmetric TISGM $\mu^*_+$ is
$\omega^{(n)}(x)\equiv 1$.
\vskip 0.5 truecm

The boundary conditions for the bottom and middle
symmetric TISGM, $\mu^*_-$ and $\mu_{\rm m}$, are unclear.
In the case of a general $m$, we also have two conjectures.
\vskip 0.5 truecm

{\bf Conjecture 5.} $\forall$ $m,k\geq 2$ and $J<0$,
there exist symmetric solutions $h=(h_0,h_1,...,h_{m-1})$ to
(2.6), with $h_0=0$ and $h_i=h_{m-i},\;i=1,2,...,m-1$.
\vskip 0.5 truecm

{\bf Conjecture 6.} $\forall$ $m,k\geq 2$ and $J>0$,
$\forall$ $\beta\geq 0$ the TISGM is unique
and is a symmetric TISGM.

\section{Periodic SGMs}

In this section we study a periodic (see Definition 4.1) solutions 
of system (3.1).

Note that (see [18]) there exists a one-to-one
correspondence between the set  $V$ of vertices of the CT 
of order $k\geq 1$ and the group $G_{k}$ of the free products of
$k+1$ cyclic  groups  of the second order with generators
$a_1,a_2,...,a_{k+1}$.
\vskip 0.5 truecm

{\bf Definition 4.1.} {\sl Let $K$ be a subgroup of $G_k$. 
We say that a collection (of functions) 
$h=\{h_x\in R^2 : x\in G_k\}$ is {\it $K$-periodic}
if $h_{yx}=h_x$ for all $x\in G_k$ and $y\in K$.}
\vskip 0.5 truecm

{\bf Definition 4.2.} {\sl A Gibbs measure is called 
{\it $K$-periodic} if it corresponds to
$K$-periodic collection $h$.}
\vskip 0.5 truecm

Observe that a TIGM is
$G_k$-periodic.

We give a complete description of periodic GMs i.e. a characterisation 
of such measures with respect to any normal subgroup of finite index
in $G_k$.

Let $K$ be a  subgroup of index $r$ in $G_k$, and let 
$G_k/{K}=\{K_0,K_1,...,K_{r-1}\}$ be the quotient group, with 
the coset $K_0=K$. 
Let $q_i(x)=|S_1(x)\cap K_i|, \ \ i=0,1,...,r-1$; 
$N(x)=|\{j:q_j(x)\ne 0\}|,$ where 
$S_1(x)=\{y\in G_k: \langle x, y \rangle \},$
 $x\in G_k$ and $|\cdot|$ is the number of elements in the set.
 Denote $Q(x)=(q_0(x),q_1(x),...,q_{r-1}(x)).$
 
We note (see [21]) that for every $x\in G_k$ there is a permutation 
$\pi_x$ of the coordinates of the vector $Q(e)$ (where $e$ is the 
identity of $G_k$) such that 
$$ \pi_xQ(e)=Q(x). \eqno(4.1)$$ 
It follows from (4.1) that $N(x)=N(e)$ for all $x\in G_k.$

Each $K-$ periodic collection is given by 
$$\{h_x=h_i \ \ \mbox{for} \ \ x\in K_i,\ \  i=0,1,...,r-1\}.$$

By Proposition 2.1 (for $m=2$) and (4.1), vector 
$h_n, \ \ n=0,1,...,r-1$, 
satisfies the system
$$h_n=\sum_{j=1}^{N(e)}q_{i_j}(e)F(h_{\pi_n(i_j)};\theta)-
F(h_{\pi_n(i_{j_0})};\theta), \eqno(4.2)$$
where $j_0=1,...,N(e)$, and 
function $h\mapsto F(h,m,\theta)$ defined in Proposition 2.1
takes now the form $h\mapsto F(h)=(F_0(h,\theta)$,
$F_1(h,\theta))$
where
$$F_0(h,\theta )=\ln{\exp(h_0)+\theta \exp(h_1)+\theta^2 \over
\theta^2\exp(h_0)+\theta \exp(h_1)+1},$$
$$F_1(h,\theta)=\ln{\theta\exp(h_0)+\exp(h_1)+\theta \over
\theta^2\exp(h_0)+\theta \exp(h_1)+1}.\eqno(4.3)$$
Recall, $\theta$ has been defined in (2.5).
\vskip 0.5 truecm

{\bf Proposition 4.3.} {\sl If $\theta\ne 1$, then $F(h)=F(l)$ if and
only if $h=l$.}
\vskip 0.5 truecm

{\bf Proof.} {\it Necessity.} From $F(h)=F(l)$ we get the system 
of equations
$$\left\{\begin{array}{ll}
\theta\Big(\exp(h_0+l_1)-\exp(h_1+l_0)\Big)+
(1+\theta^2)\Big(\exp(h_0)-\exp(l_0)\Big)+
\theta \Big(\exp(h_1)-\exp(l_1)\Big)=0,\\
\theta\Big(\exp(h_0)-\exp(l_0)\Big)+ \exp(h_1)-\exp(l_1)=0.\\
\end{array}\right.\eqno (4.4)$$
where $h=(h_0,h_1),\ \  l=(l_0,l_1).$
 Using the fact that

$$ \exp(h_0+l_1)-\exp(h_1+l_0)=
\exp(l_1)\Big(\exp(h_0)-\exp(l_0)\Big)-
\exp(l_0) \Big(\exp(h_1)-\exp(l_1)\Big),$$
we obtain 

$$\left\{\begin{array}{ll}
\Big(1+\theta^2 +\theta\exp(l_1)\Big)\Big(\exp(h_0)-\exp(l_0)\Big)+
\theta \Big(1-\exp(l_0)\Big)\Big(\exp(h_1)-\exp(l_1)\Big)=0,\\
\theta\Big(\exp(h_0)-\exp(l_0)\Big)+ \exp(h_1)-\exp(l_1)=0.\\
\end{array}\right.\eqno (4.5)$$

From (4.5) we get

$$\Big(1+\theta^2\exp(l_0)+\theta\exp(l_1)\Big)\Big(\exp(h_0)-
\exp(l_0)\Big)=0. \eqno(4.6)$$
It follows from (4.6) that $h_0=l_0.$ Consequently, from 
second equation in (4.5) we have $h_1=l_1$.
\vskip 0.5 truecm

{\it Sufficiency.} Straightforward.

\vskip 0.5 truecm

Let $G^*_k$ be the subgroup in $G_k$ consisting of all words 
of even length. Clearly, $G^*_k$ is a subgroup of index 2.
\vskip 0.5 truecm

{\bf Theorem 2.}  {\sl Let $K$ be a normal subgroup of finite 
index in $G_k$. Then each $K-$ periodic GM for SOS model
 is either TI or $G^*_k-$ periodic.}  
\vskip 0.5 truecm

{\bf Proof.} We see from (4.2) that 
$$F(h_{\pi_n(i_1)})=F(h_{\pi_n(i_2)})=...=F(h_{\pi_n(i_{N(e)})}).$$
Hence from Proposition 4.3 we have
$$h_{\pi_n(i_1)}=h_{\pi_n(i_2)}=...=h_{\pi_n(i_{N(e)})}.$$
Therefore,
$$ h_x=h_y=h, \ \ \mbox{if} \ \  x,y\in S_1(z),\ \ z\in G^*_k;$$
$$ h_x=h_y=l, \ \ \mbox{if} \ \  x,y\in S_1(z),\ \  
z\in G_k\setminus G^*_k.$$
Thus the measures are TI (if $h=l$) or $G^*_k-$ periodic (if $h\ne l$).
This completes the proof of Theorem  2. 
\vskip 0.5 truecm

Let $K$ be a  normal subgroup of finite index in $G_k.$ 
What condition on $K$ will guarantee that each $K-$periodic 
GM is TI? We put $I(K)=K\cap \{a_1,...,a_{k+1}\},$ 
where $a_i, \ \ i=1,...,k+1$ are generators of $G_k$.
\vskip 0.5 truecm

{\bf Theorem 3} {\sl If $I(K)\ne \emptyset$, then each 
$K-$ periodic GM  for SOS model is TI.}
\vskip 0.5 truecm

{\bf Proof.} Take $x\in K.$ We note that the inclusion 
$xa_i\in K$ holds if and only if $a_i\in K.$ 
Since $I(K)\ne \emptyset$, there is an element $a_i\in K.$ 
Therefore $K$ contains the subset $Ka_i=\{xa_i: x\in K\}$. 
By Theorem 2 we have $h_x=h$ and $h_{xa_i}=l.$ 
Since $x$ and $xa_i$ belong to $K$, it follows that 
$h_x=h_{xa_i}=h=l.$ Thus each $K-$ periodic GM is TI.
This proves Theorem 3.
\vskip 0.5 truecm

Theorems 2 and 3  reduce the problem of describing $K-$ periodic 
GM  with $I(K)\ne \emptyset$ to describing the fixed 
points of $kF(h;\theta)$ (see (3.2.a,b)) which describs 
TIGM. If $I(K)=\emptyset$, this problem is redused to 
describing the solutions of the system: 
$$\left\{\begin{array}{ll}
h=kF(l;\theta),\\
l=kF(h;\theta).\\
\end{array}\right.\eqno (4.7)$$

Denote $z_i=\exp(h_i),\ \  t_i=\exp(l_i), \ \ i=0,1.$ 
Then  from (4.7) we get
$$\left\{\begin{array}{llll}
z_0=\bigg({t_0+\theta t_1+\theta^2 \over \theta^2t_0+
\theta t_1+1}\bigg)^k,\\
z_1=\bigg({\theta t_0+ t_1+\theta \over \theta^2t_0+
\theta t_1+1}\bigg)^k,\\
t_0=\bigg({z_0+\theta z_1+\theta^2 \over \theta^2z_0+
\theta z_1+1}\bigg)^k,\\
t_1=\bigg({\theta z_0+ z_1+\theta \over \theta^2z_0+
\theta z_1+1}\bigg)^k.\\
\end{array}\right.\eqno (4.8)$$
\vskip 0.5 truecm

{\bf Proposition 4.4.} {\sl For a ferromagnetic SOS model, with $J<0$
$(\theta <1)$ (and even for $J=0$), the system of equations} 
(4.8) {\sl has solutions with $z_0=t_0$ 
and $z_1=t_1$ only.}
\vskip 0.5 truecm

{\bf Proof.} Denote $u_i=z_i^{1/k}, \ \ v_i=t_i^{1/k}, \ \ i=0,1.$
Then from (4.8) we have
$$u_0-v_0={(1-\theta^2)[\theta(u_1^kv_0^k-u_0^kv_1^k)+
(\theta^2+1)(v_0^k-u_0^k)+\theta(v_1^k-u_1^k)]\over 
(\theta^2v_0^k+\theta v_1^k+1)(\theta^2 u_0^k+\theta u_1^k+1)} \eqno (4.9)$$
and
$$u_1-v_1={(1-\theta^2)[\theta (v_0^k-u_0^k)+(v_1^k-u_1^k)]\over 
(\theta^2v_0^k+\theta v_1^k+1)(\theta^2 u_0^k+\theta u_1^k+1)} \eqno (4.10)$$
Using the fact that
$$u_1^kv_0^k-u_0^kv_1^k=
u^k_1(v_0^k-u_0^k)+u_0^k(u_1^k-v_1^k)$$
we obtain
$$\left\{\begin{array}{ll}
[A+(1-\theta^2)(\theta u^k_1+\theta^2+1)B_0](u_0-v_0)+
\theta (1-\theta^2)(1-u_0^k)B_1(u_1-v_1)=0,\\
\theta(1-\theta^2)B_0(u_0-v_0)+
[A+(1-\theta^2)B_1](u_1-v_1)=0,\\
\end{array}\right.\eqno (4.11)$$
where $$A=(\theta^2v_0^k+\theta v_1^k+1)(\theta^2 u_0^k+\theta u_1^k+1)>0,$$
$$B_i=u_i^{k-1}+u_i^{k-2}v_i+...+v_i^{k-1}>0, \ \ i=0,1.$$
From (4.11) we get
$$\Big[A^2+(1-\theta^2)\Big(B_1+(\theta u_1^k+\theta^2+1)B_0\Big)A+
(1-\theta^2)^2(\theta u_1^k+
\theta^2u_0^k+1)B_0B_1\Big](u_0-v_0)=0.\eqno(4.12)$$
Since $\theta\leq 1$ $(J\leq 0)$, we deduce from (4.12) 
that $u_0=v_0.$ Then from second equation of (4.11) we have $u_1=v_1$.
This completes the proof.
\vskip 0.5 truecm

Now  consider anti-ferromagnetic case, with $J>0 \ \ (\theta>1).$
By Proposition 3.2 we know that if $J>0$ then the system of equations
(4.8) has a unique solution with $z_0=t_0$, $z_1=t_1$. Moreover, $z_0=1$.
 For $z_0=t_0=1$ 
from (4.8) we have
$$\left\{\begin{array}{ll}
z_1=\bigg({2\theta + t_1 \over \theta^2+
\theta t_1+1}\bigg)^k,\\
t_1=\bigg({2\theta + z_1 \over \theta^2+
\theta z_1+1}\bigg)^k.\\
\end{array}\right.\eqno (4.13)$$
The following proposition gives a condition under which (4.8) has 
solutions with  $z_0=t_0=1$ and $z_1\ne t_1.$  

\vskip 0.5 truecm

{\bf Proposition 4.5.} {\sl Let $(z_*,z_*)$ be the unique solution of}
 (4.13). {\sl If 
$$ {kz_*(\theta^2-1)\over (2\theta+z_*)(1+\theta^2+\theta z_*)}>1,
\eqno(4.14)$$
then the system of equations} (4.13) {\sl has  at least three solutions
  $(z_-^*,z_+^*), (z_*,z_*), (z_+^*,z_-^*)$, where 
$z_-^*=\psi(z^*_+,\theta,k)$ and}
$$\psi(x,\theta,k)=\bigg({2\theta+x\over 1+\theta^2+\theta x}\bigg)^k.$$
\vskip 0.5 truecm

{\bf Proof.} Under (4.14) $z_*$ is unstable fixed point of the map 
$z>0\to \psi(z,\theta,k).$ For any $z\geq 1$, iterates 
$\psi^{(2n)}(z,\theta,k)$ remain $>z_*$ monotonically decrease
and hence converge to a limit, $z^*_+\geq z_*$ which solves 
$$z=\psi(\psi(z,\theta,k),\theta,k).\eqno (4.15)$$
However, $z^*_+> z_*$ as $z_*$ is unstable. Then 
$z^*_-=\psi( z_+^*,\theta,k)$ is $<z_*$ and also solves (4.15). 
This completes the proof.  
\vskip 0.5 truecm

Summarising, we obtain the following
\vskip 0.5 truecm
 
{\bf Theorem 4.} {\sl For the SOS model  with respect to any  normal 
subgroup $K\subset G_k$ of finite index the following assertions
hold:

(i) In the FM case ($J<0$), and for $J=0$ (no interaction), 
the $K$-periodic GMs coincide with TIGMs.

(ii) In the AFM case ($J>0$): (a) if  $I(K)\ne \emptyset$ then $K$-periodic
GMs coincide with TIGMs; (b) if} (4.14) {\sl holds and $I(K)=\emptyset$
then there are three $K$-periodic GMs
$\mu_{12}$, $\mu_{21}$ and $\mu_*$. Moreover, measure $\mu_*$ is TI and
mesures $\mu_{12}$ and $\mu_{21}$ are $G^*_k-$ periodic.}\\
\vskip 0.5 truecm

\section{Non-periodic SGMs}

In this section we consider the case $J<0$, $m=2$, $\beta>\beta_{\rm cr}^1$.
We use measures $\mu^*_i$, $i=-,{\rm m},+$, to show that system
(3.1) admits uncountably many non-periodic
solutions.

Take an arbitrary infinite path $\pi=\{x_0,x_1,...\}$ on
the CT ${\cal T}^k$ starting at the origin $x^0$:
$x_0=x^0$. We
will establish a 1-1 correspondence between such
paths and real numbers $t\in [0;{k+1\over k}]$ (cf. Ref. [16,17]).
In fact, let $\pi_1=\{x_0,x_1,...\}$ and $\pi_2=\{y_0,y_1,...\}$
be two such paths, with $x_0=y_0=x^0$. We will map the
pair $(\pi_1, \pi_2)$ to
a vector-function $h^{\pi_1\pi_2}:\;x\in V\mapsto h_x^{\pi_1\pi_2}$
satisfying (3.1). Paths $\pi_1$ and $\pi_2$ split
${\cal T}^k$ into three components  ${\cal T}_1^k$, ${\cal T}_2^k$
and ${\cal T}_3^k$
when $\pi_1$, $\pi_2$ are distinct and into two components ${\cal T}_1^k$
and ${\cal T}_3^k$ when $\pi_1$, $\pi_2$ coincide (again cf. Ref. [16,17]).
Vector-function $h^{\pi_1\pi_2}$ is then defined by
$$h_{x}^{\pi_1\pi_2}=\left\{\begin{array}{lll}
h^*_-,& \mbox{if} &x\in {\cal T}^k_1,\\
h^*_{\rm m},& \mbox{if} &x\in {\cal T}^k_2,\\
h^*_+,& \mbox{if} &x\in {\cal T}^k_3,\\
\end{array}\right.\eqno (5.1)$$
where vectors $h_i^*=(0,\ln z_{1,i}^*)$, $i=-,{\rm m},+$,
are solutions of (3.1).

Let $h=(h_0,h_1)\in R^2$. Denote
$$ \|h\|=\max\{|h_0|,|h_1|\}.$$
Let function $h\mapsto F(h)=F(h,\theta)$ be defined by (4.3).
\vskip 1 truecm

{\bf Proposition 5.1.} {\sl For any $h=(h_0,h_1)\in R^2$ the
following inequalities hold:}

a)
$$\left|{\partial F_i\over \partial h_j}\right|
\leq {|\theta^2-1|\over \theta^2},\;\;i,j=0,1,$$

b)
$$ \|F(h,\theta )-F(l,\theta )\|\leq
2{|\theta^2-1|\over \theta^2}\|h-l\|,\;\;h,l\in R^2,$$

c) {\sl for any $h=(0,h_1)$ and $l=(0,l_1)$:}
$$ \|F(h)-F(l)\|\leq {|\theta^2-1|\over 1+
3\theta^2+2\theta\sqrt{2(\theta^2+1)}}
\|h-l\|,\;\; h,l\in R^2.$$

d)
$$ |F_0(h)|\leq {|\theta^2-1|\over \theta^2+1}|h_0|,\;\;
h=(h_0,h_1)\in R^2.$$
\vskip 1 truecm

{\bf Proof.}
a) Write:
$${\partial F_0\over \partial h_0}=
{(1-\theta^2)e^{h_0}(\theta e^{h_1}+\theta^2+1)
\over (e^{h_0}+\theta e^{h_1}+\theta^2)(\theta^2 e^{h_0}+\theta e^{h_1}+1)}.$$
\\
 To assess the derivative ${\partial F_0\over \partial h_0}$,
consider two cases:

Case 1: $h_0\geq 0$. Then
$${\theta e^{h_1}+\theta^2+1
\over e^{h_0}+\theta e^{h_1}+\theta^2}\leq 1, \ \
{e^{h_0}
\over \theta^2 e^{h_0}+\theta e^{h_1}+1}<{1\over \theta^2}.$$

Case 2: $h_0\leq 0$. Then
$${ 1
\over e^{h_0}+\theta e^{h_1}+\theta^2}\leq {1\over \theta^2}, \ \
{(\theta e^{h_1}+\theta^2+1)e^{h_0}
\over \theta^2 e^{h_0}+\theta e^{h_1}+1}\leq 1.$$
Hence, $|{\partial F_0\over \partial h_0}|\leq {|1-\theta^2|\over
\theta^2}.$

To assess ${\partial F_0\over \partial h_1}$, we again
consider two cases:

Case 3: $h_0\geq 0$. Then
$$\left|{dF_0\over dh_1}\right|=|\theta^2-1|\theta\; {e^{h_1}
\over e^{h_0}+\theta e^{h_1}+\theta^2}\; {e^{h_0}-1\over
\theta^2 e^{h_0}+\theta e^{h_1}+1}\leq {|1-\theta^2|\over \theta^2}.$$

Case 4: $h_0<0$. Then
$$\left|{dF_0\over dh_1}\right|=|\theta^2-1|\theta \;{e^{h_1}
\over \theta^2e^{h_0}+\theta e^{h_1}+1}\; {1-e^{h_0}\over
e^{h_0}+\theta e^{h_1}+\theta^2}\leq {|1-\theta^2|\over \theta^2}.$$

Finally, to assess the derivatives of $F_1$,
write:
$$\left|{\partial F_1\over \partial h_0}\right|=
\left|\theta (\theta^2-1) {1
\over \theta e^{h_0}+ e^{h_1}+\theta}\; {e^{h_0}\over
\theta^2 e^{h_0}+\theta e^{h_1}+1}\right|\leq {|1-\theta^2|\over \theta^2}.$$
and
$$\left|{\partial F_1\over \partial h_1}\right|=\left|(\theta^2-1)
{1\over \theta e^{h_0}+ e^{h_1}+\theta}\; {e^{h_1}\over
\theta^2 e^{h_0}+\theta e^{h_1}+1}\right|\leq
{|1-\theta^2|\over\theta^2}.$$

b) Write:
 $$ \|F(h)-F(l)\|=\max\{|F_0(h)-F_0(l)|,|F_1(h)-F_1(l)|\}\leq$$
$$\max_{i=0,1}\{|(F_i)'_{h_0}||h_0-l_0|+|(F_i)'_{h_1}||h_1-l_1|\}\leq
 2{|\theta^2-1|\over \theta^2}\|h-l\|.$$

In cases c) and d) the inequalites are straightforward.
This completes the proof of Proposition 5.1.
\vskip 1 truecm

With the help of Proposition 5.1 it is easy to prove
the following Theorem 5, similar to Theorem 3 of [17]:
\vskip 1 truecm

{\bf Theorem 5.} {\sl For any two infinite paths $\pi_1$,
$\pi_2$, there
exists a unique vector-function $h^{\pi_1\pi_2}$ satisfying}
(3.1) {\sl and} (5.1).
\vskip 1 truecm

Next, we map $(\pi_1,\pi_2)$ to a pair
$(t,s)\in [0,{k+1\over k}]$x$[0,{k+1\over k}]$. In the
standard way (see [5, 16-18]) one can prove that functions
$h^{\pi_1(t)\pi_2(s)}$ are different
for different pairs $(t,s) \in D$ where
$D=\{(u,v)\in [0,{k+1\over k}]^2:u\leq v\}.$

Now let $\mu(t,s)$ denote the SGM corresponding
to function $h^{\pi_1(t)\pi_2(s)}$, $(t,s)\in D.$ We obtain the
following
\vskip 1 truecm

{\bf Theorem 6.} {\sl For any pair $(t,s)\in D$, there exists
a unique SGM $\mu(t,s)$. Moreover, the above
GMs $\mu^*_i$, $i=-,{\rm m},+$, are specified as}
$\mu(0,0)=\mu^*_+$, $\mu(0,{k+1\over k})=\mu^*_{\rm m}$,
$\mu({k+1\over k},{k+1\over k})=\mu^*_-.$
\vskip 1 truecm

Because measures $\mu(t,s)$ are different for different
$(t,s)\in D$ we obtain a continuum of distinct extreme 
SGMs.

Concluding the paper, we state our final (an perhaps most ambitious)
conjecture:
\vskip 1 truecm

{\bf Conjecture 7.} For a general ferromagnetic SOS model ($J<0$),
for temperature $T>0$ small enough, there exists at least
three translation-invariant SGMs for $m>1$ even and
at least four for $m>1$ odd. The precise numbers may depend
on $m$.
\vskip 1 truecm

{\bf Acknowledgements.} UAR thanks Cambridge Colleges Hospitality
Scheme for supporting the visit to Cambridge in July, 2002.
YMS worked in association with the ESF/RSDES Programme
``Phase Transitions and Fluctuation Phenomena
for Random Dynamics in Spatially Extended
Systems'' and was supported by the INTAS Grant 0265
Mathematics of Stochastic Networks.
UAR and YMS thank
IHES, Bures-sur-Yvette, and the IGS programme at
the Isaac Newton Institute, University of Cambridge,
for support and hospitality.
\vskip 1 truecm

{\bf References}
\vskip 1 truecm

1. Ya.G. Sinai, {\it Theory of phase transitions: Rigorous Results}
(Pergamon, Oxford, 1982).

2. C. Preston, {\it Gibbs states on countable sets} (Cambridge University
Press, London 1974); F. Spitzer, Markov random fields on an
infinite tree, {\it Ann. Prob.} {\bf 3}: 387--398 (1975).

3. V.A. Malyshev, R.A. Minlos. {\it Gibbs random fields} (Nauka, Moskow 1985).

4. H.O. Georgii, {\it Gibbs measures and phase transitions} (Walter de
Gruyter, Berlin, 1988).

5. S. Zachary, Countable state space Markov random fields and Markov
chains on trees. {\it Ann. Prob.} {\bf 11}: 894--903 (1983).

6. S. Zachary, Bounded, attractive and repulsive Markov
specifications on trees and on the one-dimensional lattice.
{\it Stochastic Process. Appl.} {\bf 20}:247--256 (1985).

7. P.M. Bleher, N.N. Ganikhodjaev, On pure phases of the Ising model
on the Bethe lattice, {\it Theor. Probab. Appl.} {\bf 35}: 216-227 (1990).

8. P.M. Bleher, Extremity of the disordered phase in the Ising model on the
Bethe lattice, {\it Comm. Math. Phys.} {\bf 128}: 411-419 (1990).

9. P.M. Bleher, J. Ruiz, V.A. Zagrebnov, On the purity of the limiting
Gibbs state for the Ising model on the Bethe lattice, {\it Journ. Statist.
Phys}. {\bf 79}:473-482 (1995).

10. P.M. Bleher, J. Ruiz, V.A. Zagrebnov, On the phase diagram of the
random field Ising model on the Bethe lattice, {\it Journ. Statist. Phys}.
{\bf 93}: 33-78 (1998).

11. D. Ioffe, On the extremality of the disordered state for
the Ising model on the Bethe lattice, {\it  Lett. Math. Phys}.
{\bf 37}: 137-143 (1996).

12. D. Ioffe.  Extremality of the disordered state for the Ising
model on general trees, { \it Trees}, Versailles, 1995. {\it Progr. Probab.}
{\bf 40}:3-14 (Progr. Probab., 40, Birkhauser, Basel, 1996).

13. P.M. Bleher, J. Ruiz, R.H. Schonmann, S. Shlosman, V.A. Zagrebnov,
Rigidity of
the critical phases on a Cayley tree, {\it Moscow Mathematical Journ.}.
{\bf 3}: 345-363 (2001).

14. A.E. Mazel, Yu.M. Suhov. Random surfaces with two-sided
constraints: an application of the theory of dominant ground states,
{\it Journ. Statist. Phys.} {\bf 64}:111-134 (1991).

15. F.P. Kelly, Stochastic models of computer communication systems.
With discussion. {\it Journ. Roy. Statist. Soc.} Ser B. {\bf 47}: 379--395,
415-- 428 (1985); K. Ramanan, A. Sengupta, I. Ziedins and P.
Mitra, Markov random field models of multicasting in
tree networks, {\it Adv. Appl. Probab.} {\bf 34}:58--84 (2002).

16. U.A. Rozikov,  Describtion of limiting Gibbs measures
for $\lambda-$ models on the Bethe lattice, {\it Siberian Math.
Journ.} {\bf 39}:427-435 (1998).

17. U.A. Rozikov,  Describtion uncountable number of Gibbs
measures for inhomogeneous Ising model, {\it Theor. Math. Phys.}
{\bf 118}:95-104 (1999).

18. N.N. Ganikhodjaev, U.A. Rozikov, Describtion of periodic
extreme Gibbs measures of some lattice models on the Cayley
tree, {\it Theor. Math. Phys}.{\bf 111}: 480-486 (1997).

19. N.N. Ganikhodjaev, U.A. Rozikov, On disordered phase in the ferromagnetic
Potts model on the Bethe lattice, {\it Osaka Journ. Math.}.{\bf 37}:373-383
(2000).

20. N.N. Ganikhodjaev, On pure phases of the ferromagnet Potts with
three states on the Bethe lattice of order two. {\it Theor. Math.
Phys.} {\bf 85}:163--175 (1990).

21. U.A.Rozikov, Partition structures of the group representation of the
Cayley tree into cosets by finite-index normal subgroups and their 
applications to the description of periodic Gibbs distributions.
{\it Theor. Math. Phys.} {\bf 112}: 929-933 (1997).
\end{document}